\documentclass[11pt]{amsart}

\usepackage{graphicx}
\usepackage{amssymb}

\usepackage{amsmath}

\newcommand{\x}{\dot{x}}

\newenvironment{pf*}[1]{\begin{proof}}{\end{proof}}
\newtheorem{thm}{Theorem}
\newtheorem{lem}[thm]{Lemma}
\newtheorem{cor}[thm]{Corollary}
\newtheorem{prop}[thm]{Proposition}

\newcommand{\idiot}[1]{\vspace{5 mm}\par \noindent
\marginpar{\textsc{Note}}
\framebox{\begin{minipage}[c]{0.95 \textwidth}
#1 \end{minipage}}\vspace{5 mm}\par}

\renewcommand{\idiot}[1]{}

\newcommand{\shortorlongversion}[2]{#2}

\def\la{{\lambda}}
\def\shuf{{\sqcup\!\sqcup}}
\def\QQ{{\mathbb Q}}

\def\CC{{\mathbb Q}}
\def\lab{{\langle}}
\def\rab{{\rangle}}
\def\S{{\mathfrak S}}
\def\frobn{{{\mathcal Frob}_{\S_n}}}
\def\frobd{{{\mathcal Frob}_{\S_d}}}
\def\LL{{\mathcal L}}
\def\HH{{\mathcal H}}
\def\AA{{\mathcal A}}

\newcommand{\Sym}{\textit{Sym}}
\newcommand{\pchoose}[2]{\begin{pmatrix}#1\\ #2\end{pmatrix}}

\begin{document}


\title[Noncommutative harmonics]{ On the  $\S_n$-module structure of the noncommutative harmonics. }


\author{Emmanuel Briand}
\address{Emmanuel Briand and Mercedes Rosas, Universidad de Sevilla, Sevilla, Spain}

\author{Mercedes Rosas}

\author{Mike Zabrocki}
\address{Mike Zabrocki, York University, Toronto, Canada}

\thanks{Expanded version of a paper to appear in \emph{Journal of Combinatorial theory, series A}, \texttt{http://www.elsevier.com/locate/jcta}.\\
\indent Emmanuel Briand is supported by a contract \emph{Juan de la Cierva}, MEC. Mercedes Rosas is supported by a contract \emph{Ram{\'o}n y Cajal}, MEC. Mike Zabrocki is supported by NSERC}


\begin{abstract}
Using a noncommutative analog of Chevalley's decomposition of polynomials into symmetric polynomials times coinvariants due to 
Bergeron, Reutenauer, Rosas, and Zabrocki we compute the
graded Frobenius characteristic for their two sets of noncommutative harmonics
with respect to the left action of the symmetric group (acting on variables). 
We  use these results to derive the Frobenius
series for the enveloping algebra of the derived free Lie algebra
in $n$ variables. 
\end{abstract}



\maketitle

{\em In honor of Manfred Schocker (1970-2006).  The authors would also
like to acknowledge the contributions that he made to this paper.}

\section{Introduction}





A central result of Claude Chevalley \cite{C}
decomposes the ring of polynomials in $n$ variables (as graded representation of the symmetric group $\S_n$)
 as the tensor product of the symmetric polynomials 
times the coinvariants of $\S_n$
 (i.e., polynomials modulo symmetric polynomials with no constant term). 
 
The coinvariants of the symmetric group can also be defined as its harmonics (the polynomials annhilated by all symmetric polynomial differential operators with no constant term). 
They admit as a basis the famous
Schubert polynomials of Schubert  calculus, 
that  play an 
important role in algebraic combinatorics, see for instance \cite{F}. 
 
The space of invariant polynomials in  noncommutative variables 
was introduced in 1936 by Wolf \cite{W}
where she found a noncommutative version of the fundamental
theorem of symmetric functions. This space has been 
studied from a modern perspective in \cite{RS, BRRZ, BZ}. On the other hand,
two sets of noncommutative harmonics for the symmetric group 
were introduced in \cite{BRRZ} that translated into two 
noncommutative  analogues of Chevalley decomposition for the
ring of polynomials in noncommuting variables.  
The question of decomposing as
$\S_n$--modules both kinds of noncommutative harmonics 
 was left open. This is the starting point in our investigations. 

We begin the present work with the computation of
the graded Frobenius characteristic
of noncommutative harmonics.
We then use these calculations to derive the Frobenius
series for the enveloping algebra of the derived free Lie algebra
in $n$ variables, $\AA_n'$. 
This last computation is achieved by using the existence of an isomorphism of 
$GL_n(\QQ)$--modules 
between the space of polynomials
in noncommutative variables, and the tensor product of the
space of commuting polynomials with 
$\AA'_n$. \shortorlongversion{}{Such an isomorphism is presented explicitly in the last section.}

We conclude this introduction with some basic definitions and 
results that we will be using in the following sections. Let $\S_n$ denote the symmetric group in $n$ letters.
Denote by $\QQ[X_n] = \QQ[x_1, x_2, \ldots, x_n]$  the space of
polynomials in $n$ commuting variables
  and by  $\QQ\lab X_n\rab = \QQ\lab x_1, x_2, \ldots, x_n\rab$ 
the space of polynomials in $n$ noncommutative
 variables.

The space of symmetric polynomials in $n$ variables will be denoted
by $Sym_n$ and the space of noncommutative polynomials which are invariant under
the canonical action of the symmetric group $\S_n$ will be denoted by
$NCSym_n$.

Given any  polynomial $f(X_n) \in \QQ[X_n]$, the notation $f(\partial_{X_n})$
represents the polynomial turned into an operator with each of the
variables replaced by its corresponding derivative operator.  Analogous
notation will  also hold for $f(X_n) \in \QQ\lab X_n \rab$ except that 
there are two types of differential operators acting on words in noncommutative
variables.  The first is the {\it Hausdorff derivative}, $\partial_x$, whose action on a
word $w$ is defined to be the sum of the subwords of $w$ with an occurrence
of the letter $x$ deleted.  The second derivative is the {\it twisted derivative}, $d_x$,
which is defined on $w$ to be  $w'$ if $w = xw'$, and $0$
otherwise. Both derivations are extended to polynomials by
linearity. 

It is interesting to remark (as does Lenormand in \cite{L}, section \emph{S{\'e}ries comme op{\'e}rateurs}) that these two
operations are dual to the shuffle and concatenation products respectively,
with respect to a scalar product where the noncommutative monomials are
self dual.  That is,
\begin{eqnarray*}
\lab \partial_x f, g \rab &=& \lab f, x \shuf g \rab, \hbox{ and } \\
\lab d_x f, g \rab &=& \lab f, x g \rab.
\end{eqnarray*}


Following \cite{BRRZ}, we introduce the following two sets of
 noncommutative analogues of the harmonic polynomials. The canonical
action of the symmetric group endow them with the structure of
$\S_n$--modules. 
\begin{eqnarray*}
MHar_n\, &=& \{ f \in \QQ\lab X_n\rab :\; p(\partial_{X_n}) f(X_n) = 0\hbox{ for all }
p \in \mathfrak{M}_n \}\\
NCHar_n &=& \{ f \in \QQ\lab X_n\rab :\; p(d_{X_n}) f(X_n) = 0\hbox{ for all }
p \in \mathfrak{M}_n \}
\end{eqnarray*}
where
$
\mathfrak{M}_n=\{
p \in NCSym_n\hbox{ with }p(0)= 0\}
$.

We are now ready to state the two decompositions of  $\QQ\lab X_n
\rab$ as the tensor product (over $\QQ$) of its invariants times its coinvariants
that we have described.
\begin{prop}[\cite{BRRZ}, Theorems 6.8 and 8.8] \label{harmoniciso}\label{hariso} 
As  graded $\S_n$--modules,
\begin{eqnarray*}
 \QQ \lab X_n \rab &\simeq& \,\, MHar_n \,\, \otimes \,\, Sym_n, \\
 \QQ \lab X_n \rab &\simeq& NCHar_n \otimes NCSym_n. 
\end{eqnarray*}
\end{prop}

\begin{section}{ The Frobenius characteristic of noncommutative harmonics}
In this section we compute the Frobenius characteristic of  both kinds of
 noncommutative harmonics.
This section is based of the observation that the graded Frobenius
series for each of the $\S_n$--modules appearing in
Proposition \ref{harmoniciso} is either known or 
can be deduced from the existence of the isomorphisms described
there. 

The expressions for Frobenius images and characters will require a little
use of symmetric function notation and identities. 
We will follow Macdonald \cite{M} for the notation
of the $s_\la$ Schur, $h_\la$ homogeneous, $e_\la$ elementary and
$p_\la$ power sums bases for the ring of symmetric functions $\Sym$, that
we identify with $\QQ[p_1, p_2, p_3, \ldots]$. For convenience we will make use of
some plethystic notation. 

For a symmetric function $f$, $f[X]$ represents
the symmetric function evaluated at an unspecified (possibly infinite) alphabet
$X$. Then, $f[X(1-q)]$ is the image of $f$ under the algebra automorphism sending 
the power sum symmetric function $p_k$
to $(1-q^k) p_k[X]$. Similarly,
$f\left[ \frac{X}{1-q} \right]$ is the image of the symmetric function $f$ under the inverse automorphism (sending the power sum $p_k$ to 
$p_k/(1-q^{k})$).

 In our calculations, we use the Kronecker product $\odot$ of symmetric
 functions. This operation on symmetric functions corresponds, under 
the Frobenius map, to the inner tensor product of representations of the 
symmetric group (tensor product of representations with the diagonal action on the tensors).
It can also be defined directly on symmetric
 functions by the equation $ p_{\lambda} \odot  p_{\mu} =
 \delta_{\lambda, \mu} \left(\prod_{i} n_i(\la)! i^{n_i(\la)}\right) p_{\lambda}$ where
 $n_i(\la)$ is the number of parts of size $i$ in $\la$, and then extended by
bilinearity.

We introduce the notations
\begin{eqnarray*}
(q;q)_k  &=& (1-q)(1-q^2) \cdots (1-q^k), \\
\{q;q\}_k &=& (1-q)(1-2q) \cdots (1-k\,q). 
\end{eqnarray*}
 Then $q^d/\{ q;q\}_d $ is the generating function for the set partitions
 with length $d$ and $q^d/(q;q)_d$ is the generating function for partitions
 with length $d$, \cite{Wilf}. Finally, since $Sym_n$ and $NCSym_n$ are made of graded
 copies  of the trivial $\S_n$-module we conclude that
\begin{eqnarray*}
\frobn( NCSym_n)  &=& h_n\left[ X \right] \sum_{d=0}^{n} \frac{q^d  }{ \{ q,q\}_d }\\
\frobn( Sym_n )  = h_n\left[ X \right] \frac{1}{(q;q)_n} &=& 
h_n\left[ X \right] \sum_{d=0}^{n} \frac{q^d  }{ (q;q)_d }.
\end{eqnarray*}

In the following lemma we compute the graded Frobenius characteristic for the module $\QQ\lab X_n \rab$.
\begin{lem}[The Frobenius characteristic of $\QQ\lab X_n \rab$]  \label{QXnfrob}

\noindent
\[
\frobn( \QQ\lab X_n \rab )
= \sum_{d=0}^n \frac{q^d}{ \{ q,q\}_d }h_{(n-d,1^d)}[X].
\]
\end{lem}
\idiot{
\begin{proof} (first proof)
The character of the $\S_n$--module $\QQ\lab X_n \rab$ is the trace of the action of
$\S_n$ acting on a basis and so we consider the basis of monomials in the alphabet
$X_n$.  Hence we calculate that
$\chi_q^{\QQ\lab X_n \rab}(\sigma) = \sum_{w \in X_n^*} \delta_{\sigma(w) = w} q^{|w|}$.
Since a word is fixed by the action of a permutation of $\sigma$ if and only if it
is in the alphabet of $fix(\sigma) = \#\{ i : \sigma(i) = i \}$.  
Therefore $\chi_q^{\QQ\lab X_n \rab}(\sigma) = \frac{1}{1-q fix(\sigma)}$.

Now consider the coefficient of $p_\la[X]/z_\la$ in the right hand side of \eqref{QXnfrob}.
This is calculated to be $\sum_{d=0}^n \frac{q^d}{\{q;q\}_d} d! \pchoose{fix(\lambda)}{d}$
where $fix(\la)$ is the number of parts of size $1$ in the partition $\lambda$.  A combinatorial
argument says that every word in an alphabet of size $fix(\la)$ contains exactly $d$ letters
for some $0 \leq d \leq fix(\lambda)$ (chosen
in $\pchoose{fix(\lambda)}{d}$ ways),
which appear in some order in the word (ordered in $d!$ ways)
and the placement of the positions are given by a set partition of length $d$
(the generating function for which is $\frac{q^d}{\{q;q\}_d}$).  Therefore
\begin{equation}
\sum_{d=0}^n \frac{q^d}{ \{ q,q\}_d }h_{(n-d,1^d)}[X] =
\sum_{\lambda \vdash n} \frac{1}{1-fix(\la) q} \frac{p_\la}{z_\la} = \sum_{\la \vdash n}
\chi_q^{\QQ\lab X_n \rab}(\sigma) \frac{p_\la}{z_\la} = \frobn( \QQ\lab X_n \rab ).
\end{equation}
\end{proof}
}

\begin{pf*}{Proof} 
For each monomial $x_{i_1} \cdots x_{i_r}$, we define its {\em type}
 $\nabla( x_{i_1} \cdots x_{i_r} )$ to be
the set partition of $[r] = \{1,2, \ldots, r\}$
such that $a$ and $b$ are in the same part of the set partition
if and only if $i_a = i_b$ in the monomial.  For a set partition $A$
 with at most  $n$ parts, we will let
$M^{A}$ equal the $\S_n$ submodule of $\QQ \lab X_n \rab$ spanned by all monomials
of type $A$.  As  $\S_n$--module,
\[
\QQ\lab X_n \rab \simeq \bigoplus_{d=0}^n \bigoplus_{A\, :\, \ell(A) = d}
M^A
\]
 where the second direct sum is taken over all set partitions $A$ with
 $d$ parts.

Fix a set partition $A$, and let $d$ be the number of parts of $A$, 
and ${\bf x}_{\vec i} = x_{i_1} x_{i_2} ... x_{i_r}$ be
the smallest monomial in lex order in $M^A$.  It involves only the variables $x_1$, $x_2$, \ldots, $x_d$. The representation $M^A$ is the representation of $\S_n$ induced by the action of the subgroup $\S_{d} \times \S_1^{n-d} \simeq \S_d$ on the subspace $\QQ[\S_d] \cdot {\bf x}_{\vec i}$. The representation $\QQ[\S_d] \cdot {\bf x}_{\vec i}$ of $\S_d$ is isomorphic to the regular representation.  We use the
rule for a representation $R$ of $\S_d$ induced to $\S_n$,
$$\frobn( R \uparrow^{\S_n}_{\S_d} ) = h_{n-d}[X] \frobd( R ),$$ 
and conclude that the Frobenius characteristic of $M^A$ is $h_{(n-d,1^d)}[X].$
Hence the graded Frobenius characteristic of $\QQ \lab X_n \rab$ is
\[ \frobn( \QQ \lab X_n \rab ) = \sum_{d = 0}^n \sum_{A : \ell(A) = d} q^{|A|} h_{(n-d,1^d)}[X]
= \sum_{d=0}^{n} \frac{q^d  }{ \{ q,q\}_d} h_{(n-d,1^{d})}[X]. \]
\end{pf*}

We are now able to compute the Frobenius characteristic for $MHar_n$ and
$NCHar_n$.
\begin{thm}[The Frobenius characteristic of the noncommutative harmonics] \label{frobchar}
\[ \frobn( MHar_n ) =  (q; q)_n \sum_{d=0}^{n} 
  \frac{q^d}{\{ q,q\}_d} h_{(n-d,1^{d})}[X] \]
and
\[ \frobn( NCHar_n ) = \Big (\sum_{d=0}^{n} \frac{q^d  }{\{ q,q\}_d} \Big)^{-1}
\sum_{d=0}^{n} \frac{q^d  }{ \{ q,q\}_d} h_{(n-d,1^{d})}[X]. \]
\end{thm}
\begin{pf*}{Proof}
This follows since $\frobn(MHar_n \otimes Sym_n) 
= \frobn(MHar_n) \odot \frobn(Sym_n)$.
Since $h_n[X]$ is the unity for the Kronecker product on symmetric
functions of degree $n$, and since $\frobn(Sym_n)=h_n[X]/(q;q)_n$, we
conclude that
$\frobn(MHar_n)/(q;q)_n = \frobn(\QQ\lab X_n \rab)$. We can now solve for
$\frobn(MHar_n)$.  

A similar argument demonstrates the formula for
$\frobn(NCHar_n)$.  We have from Proposition \ref{hariso} and Lemma \ref{QXnfrob},
\begin{eqnarray*}
\sum_{d=0}^n \frac{q^d}{ \{ q,q\}_d }h_{(n-d,1^d)}[X]
&=& \frobn( \QQ\lab X_n \rab )\\
&=& \frobn( NCHar_n) \odot \frobn( NCSym_n )\\
&=& \sum_{d=0}^n \frac{q^d}{ \{ q,q\}_d } h_n[X] \odot  \frobn( NCHar_n)\\
&=& \left( \sum_{d=0}^n \frac{q^d}{ \{ q,q\}_d }\right) \frobn( NCHar_n).
\end{eqnarray*}
{}From this equation we can solve for $\frobn( NCHar_n)$.
\end{pf*}

As a corollary, we obtain the generating functions for
the graded dimensions of these spaces.

\begin{cor}[The Hilbert series of the noncommutative harmonics]
\begin{eqnarray*}
dim_q( MHar_n ) &=& \frac{(q; q)_n}{1-nq}\\
dim_q(NCHar_n) &=& \frac{
1} {(1-nq)\sum_{d=0}^{n} \frac{q^d  }{\{ q,q\}_d}}
\end{eqnarray*}
\end{cor}
\begin{pf*}{Proof}
After Theorem \ref{frobchar}, 
\begin{eqnarray*}
\frobn( MHar_n ) &=&  (q; q)_n \;\frobn(\QQ\lab X_n\rab)\\
\frobn( NCHar_n ) &=& \Big (\sum_{d=0}^{n} \frac{q^d  }{\{ q,q\}_d} \Big)^{-1}
\frobn(\QQ\lab X_n\rab)
\end{eqnarray*}
This implies
\begin{eqnarray*}
dim_q( MHar_n )  &=&  (q; q)_n \;dim_q(\QQ\lab X_n\rab)\\
dim_q( NCHar_n ) &=& \Big (\sum_{d=0}^{n} \frac{q^d  }{\{ q,q\}_d} \Big)^{-1}
dim_q(\QQ\lab X_n\rab)
\end{eqnarray*}
since 
the Hilbert series of a graded $\S_n$--module is obtained by coefficient extraction from the graded Frobenius characteristic (of the coefficient of $p_{(1^n)}[X]/n!\,$ in the expansion in power sum symmetric functions).
Last, the Hilbert series of $\QQ\lab X_n\rab$ is $\frac{1}{1-nq}$.
\end{pf*}

The graded dimensions of $MHar_n$ for $2 \leq n \leq 5$ 
are listed in \cite{OEIS} as sequences
$A122391$ through $A122394$.
The sequences of graded dimensions of $NCHar_n$ 
for $3 \leq n \leq 8$ are listed in \cite{OEIS} as
sequences $A122367$  through $A122372$.

\end{section}


\begin{section}{Non--commutative harmonics and the enveloping algebra of the derived free Lie algebra}

Let  $\LL_n$  be the canonical realization of the free Lie algebra 
inside the ring of polynomials in noncommuting variables $ \QQ \lab X_n \rab$. 
More precisely,  $\LL_n$ is the linear span of the minimal set of 
polynomials in $\QQ \lab X_n \rab$ that includes $\QQ$ and the variables $X_n$, and is closed 
under the bracket operation $[x,y]=xy-yx$. Let 
$\LL'_n= [\LL_n, \LL_n]$ be the derived free Lie algebra. Remark 
that 
$\LL_n = \LL'_n \oplus \QQ X_n$, 
where $\QQ X_n$ denotes the space of linear polynomials. The
enveloping algebra $\AA_n'$ of $\LL_n$ can be realized as a subalgebra of $\QQ \lab X_n \rab$ as follows (see \cite{R} 1.6.5):
\[
\AA_n' = \bigcap_{x\in X_n} \ker \partial_x.
\]
More explicitly, $\AA_n' $ is the subalgebra of $\QQ \lab X_n \rab$  generated by all the
brackets under concatenation.
\idiot{You also said that this is in Drensky's paper.  I have found it in
`Polynomial Identities for 2x2 Matrices'  due to lemma 2.1 which we don't quote as a reference.
Perhaps it is in one of the two other references that we do quote?  I don't have
reference [D1].  My solution at the moment: add reference [D3] and quote it also for
this result.}

In \cite{BRRZ} it was established that there is an isomorphism of
vector spaces between $MHar_n$ and $ \AA_n' \otimes \HH_n$. 
In this section we will show the following result.

\begin{thm}\label{Sn iso MHar} As $\S_n$--modules,
\[ MHar_n \simeq \AA_n' \otimes \HH_n. \]
\end{thm}
The proposition will be established by comparing the Frobenius image of $MHar_n$ (known from Theorem \ref{frobchar}) to $\frobn( \AA_n' \otimes \HH_n )$, which is equal to $\frobn( \AA_n' ) \odot \frobn( \HH_n)$. We will determine $\frobn( \AA_n' )$ in Theorem \ref{frob A'n} below. An intermediate step will make use the following Theorem due to V. Drensky.

\begin{prop}[Drensky, \cite{D} Theorem 2.6] \label{glniso}
As $GL_n(\QQ)$--modules (and consequently as $\S_n$--modules),
\[ \QQ \lab X_n \rab \simeq \QQ[X_n] \otimes \AA_n'. \]
\end{prop}

\idiot{We  provide two different proofs for this theorem. For the
first proof we compute $GL_n(\QQ)$--character of
 both sides of the isomorphism and show that they 
agree , see for instance the notes by Kraft and Procesi, \cite{KP}. A second proof is provided in the final section. 

There 
 we construct an explicit
isomorphism between $\QQ \lab X_n \rab$ and $ \QQ[X_n] \otimes
\AA_n'$.
For this construction we introduce an interesting new basis for $\QQ \lab X_n \rab$
that we call the hybrid basis.}
Drensky proved Proposition \ref{glniso} by exhibiting an explicit isomorphism between these two representations. 
\shortorlongversion{
We will provide a non--constructive (but shorter) proof of the theorem.
}{
We will look at it in the next section. For now, we will provide a non--constructive proof of the theorem.
}
Before, we need to introduce some notation.

It is known that $\QQ \lab X_n \rab$ is the universal enveloping algebra
(u.e.a) 
of the free Lie algebra, $\LL_n$.
Using the Poincar\'e-Birkhoff-Witt theorem, a linear basis for $\QQ
\lab X_n \rab$
 is given by decreasing products of elements of $\LL_n$.  Since
we can choose an ordering of the elements of $\LL_n$ so that the space of linear polynomials is smallest and decreasing products of linear polynomials are
isomorphic to $\QQ [ X_n ]$ (as a vector space), we note that as vector spaces
\[ \QQ \lab X_n \rab = u.e.a.( \LL_n )  = u.e.a( \QQ X_n \oplus \LL_n' ) \simeq \QQ[X_n] \otimes \AA_n'. \]

To distinguish between the commutative elements of $\QQ [ X_n ]$
and the noncommutative words of  $\QQ \lab X_n \rab$, 
we will place a dot over the variables (as in $\x_i$)
to indicate the commutative variables.

Let  $[n] = \{ 1, 2, \ldots, n\}$
and let  $[n]^r$ denote the words of length $r$ in the alphabet
of the numbers $1, 2, \ldots, n$.  A word $w \in [n]^r$ is called a Lyndon
word if $w < w_{k} w_{k+1} \cdots w_{r}$
for all $2 \leq k \leq r$ where $<$ represents lexicographic order
on words.

Every word $w \in [n]^r$ is equal to a unique product
$w = \ell_1 \ell_2 \cdots \ell_k$ such that $\ell_1 \geq \ell_2 \geq \cdots
\geq \ell_k$ and each $\ell_i$ is Lyndon (e.g. Corollary 4.4 of \cite{R}).

Let $\ell$ be a Lyndon word of length greater than $1$.
We say that $\ell = uv$ is
the standard factorization of $\ell$ if $v$ is the smallest nontrivial
suffix in lexicographic order.  It follows that $u$ and $v$ are Lyndon
words and $u < v$.

For a Lyndon word $\ell$, if $\ell$ is a single letter $a$ then define
$P_a  = x_a \in \QQ \lab X_n \rab$.  If $\ell = uv$ is the standard factorization
of $\ell$, then $P_{\ell} = [P_u, P_v]$.  For any $w \in [n]^r$ with Lyndon
decomposition $w = \ell_1 \ell_2 \cdots \ell_k$, define
\[ P_w = P_{\ell_1} P_{\ell_2} \cdots P_{\ell_k}. \]

The set $\{ P_w \}_{w \in [n]^r}$ forms a basis
for the noncommutative polynomials of degree $r$ (\cite{R}, Theorem 5.1). The elements $P_w$ with Lyndon decomposition $w = \ell_1 \ell_2 \cdots \ell_k$ such that each Lyndon factor has degree at least $2$ are a basis of $\AA'_n$.

\begin{pf*}{Proof (of Proposition \ref{glniso})}
To prove that $\QQ \lab X_n \rab$ and $ \QQ[X_n] \otimes
\AA_n'$ are isomorphic as $GL_n(\QQ)$--modules, we
use the fact that two polynomial $GL_n(\QQ)$--modules with the same character are
isomorphic (see for instance the notes by Kraft and Procesi, \cite{KP}).  The character of a $GL_n(\QQ)$--module is the trace of
the action of the diagonal matrix $diag(a_1, a_2, \ldots, a_n)$.

A basis for $\QQ[X_n] \otimes \AA_n'$ are the elements ${\bf \x}^\alpha \otimes
P_{\ell_1} \cdots P_{\ell_k}$ with $\ell_1 \geq \ell_2 \geq \cdots \geq \ell_k$ and
$|\ell_i| \geq 2$. 
 The action of the diagonal matrix $diag(a_1, a_2, \ldots, a_n)$
on this basis element is the same as the action
on the noncommutative polynomial 
$x_1^{\alpha_1} x_2^{\alpha_2} \cdots x_n^{\alpha_n} P_{\ell_1} P_{\ell_2} \cdots P_{\ell_k}$ (in both cases: multiplication by $a_1^{\alpha_1+m_1} a_2^{\alpha_2+m_2} \cdots a_n^{\alpha_n+m_n}$ where $m_i$ is the number of occurrences of $i$ in the word $\ell_1 \ell_2 \cdots \ell_k$).
By the Poincar\'e-Birkhoff-Witt theorem, these polynomials form a basis for $\QQ\left< X_n \right>$,
hence the trace of the action of $diag(a_1, a_2, \ldots, a_n)$ acting on
$\QQ\left< X_n \right>$ and $\QQ[X_n] \otimes \AA_n'$ are equal. Since
their characters are equal, we conclude that they are isomorphic as $GL_n(\QQ)$
modules.
\end{pf*}

The $GL_n(\QQ)$--character of
$\QQ[X_n]$ is $\prod_{i=1}^n \frac{1}{1-a_i}$, and the
$GL_n(\QQ)$--character 
of $\QQ\lab X_n \rab$ is $\frac{1}{1- (a_1 + a_2 + \cdots +
  a_n)}$. Therefore,  the existence of  a 
$GL_n(\QQ)$-module isomorphism between  $\QQ \lab X_n \rab$ and $ \QQ[X_n] \otimes
\AA_n'$  implies the following result.

\begin{cor}[The $GL_n(\QQ)$--character of $\AA_n'$]
\begin{eqnarray*}
char_{GL_n(\QQ)}( \AA_n' )(a_1, a_2, \cdots, a_n) 
&=& \frac{(1-a_1) \cdots (1-a_n)}{1- (a_1 + a_2 + \cdots + a_n)}\\
&=& \sum_{k \geq 0} \sum_{i=2}^k (-1)^i e_{(i,1^{k-i})}(a_1,a_2, \ldots, a_n).
\end{eqnarray*}
Moreover this last sum is equal to
\[ \sum_T s_{shape(T)}(a_1,a_2, \ldots, a_n) \]
where the sum is over all standard tableaux $T$ such that the smallest integer
which does not appear in the first column of $T$ is odd.
\end{cor}

%
%
%
%
%
%

\idiot{
\begin{proof}
\begin{eqnarray*}
\frac{(1-a_1) \cdots (1-a_n)}{1- (a_1 + a_2 + \cdots + a_n)} &=
(\sum_{r=0}^n (-1)^r e_{r}(a_1, a_2, \ldots, a_n))(\sum_{k\geq 0} e_{(1^k)}(a_1, a_2, \ldots, a_n))\\
&= \sum_{k \geq 0} \sum_{i=0}^k (-1)^i e_{(i,1^{k-i})}(a_1,a_2, \ldots, a_n)\\
&= \sum_{k \geq 0} \sum_{i=2}^k (-1)^i e_{(i,1^{k-i})}(a_1,a_2, \ldots, a_n).
\end{eqnarray*}
\end{proof}}

By Schur-Weyl duality, the above formula also describes the
decomposition of the subspace of multilinear polynomials 
(i.e. with distinct occurrences of the variables) of
${\mathcal A}_n'$.  That is, if $n$ is the number of variables, the
the multilinear polynomials of degree $n$ will be an ${\mathfrak S}_n$-module
with Frobenius image equal to $\sum_{i=2}^n (-1)^i e_{(i,1^{n-i})}[X]$.
This decomposition was considered in the papers
\cite{D82}, \cite{P1}, \cite{P2} where an expression was given degree by degree
up to $n=7$. The expansion of this formula in the Schur basis provided in the Theorem agrees with the computations in those papers.

\idiot{In \cite{D} corollary 2.7 (i) appears explicitly the relation $\chi(\AA'_n)=\chi(\QQ\lab X_n \rab) \prod_i(1-a_i)$, more precisely the similar expressions for quotients. Remark 2.8 makes explicit the connection with proper multilinear polynomials. We think that at that moment people working on PI algebras knew enough representation theory to stop doing computations similar to those of  \cite{P2}.}

We can derive a formula for the Frobenius characteristic of
${\mathcal A}'$ by using a similar technique.

\begin{thm}[The Frobenius characteristic of $\AA_n'$]\label{frob A'n}

\noindent
\[ \frobn({\mathcal A}'_n) = 
\sum_{d=0}^{n} \frac{q^d  }{\{ q;q \}_d} h_{(n-d,1^{d})}[X(1-q)]. \]
\end{thm}
\begin{pf*}{Proof}
For any symmetric function $f[X]$ of degree $n$, we have that
\[
f[X] \odot h_n\left[ \frac{X}{1-q} \right] =
f\left[ \frac{X}{1-q}\right].
\]
  In particular, since $\frobn( \QQ[ X_n ] ) = h_n\left[ \frac{X}{1-q}
  \right],$ we conclude that
\begin{eqnarray*}
\frobn( \QQ \lab X_n \rab ) &=& \frobn( {\mathcal A}'_n \otimes \QQ[X_n])\\ 
&=& \frobn({\mathcal A}'_n) \odot h_n\left[\frac{X}{1-q}\right] 
= \frobn({\mathcal A}'_n)\left[\frac{X}{1-q}\right].
\end{eqnarray*}
This implies that if we make the plethystic substitution $X {\rightarrow} X(1-q)$ into both
sides of this equation and using Lemma \ref{QXnfrob} we arrive
at the stated formula.
\end{pf*}

We can now prove Theorem \ref{Sn iso MHar}. 
\begin{pf*}{Proof (of Theorem \ref{Sn iso MHar})}
{}From Theorem \ref{frobchar} we know the Frobenius image of $MHar_n$, we compare
this to 
\begin{eqnarray*}
\frobn( \AA_n' \otimes \HH_n ) &=& \frobn( \AA_n' ) \odot \frobn( \HH_n), \\
&=& \sum_{d=0}^{n} \frac{q^d  }{\{ q;q \}_d} h_{(n-d,1^{d})}[X(1-q)] 
\odot h_n \left[ \frac{X}{1-q} \right] (q;q)_n\\
&=& (q;q)_n \sum_{d=0}^{n} \frac{q^d  }{\{ q;q \}_d} h_{(n-d,1^{d})}[X]\\
&=& \frobn( MHar_n).
\end{eqnarray*}
Since the two $\S_n$--modules have the same Frobenius image, we conclude that they
must be isomorphic.
\end{pf*}
\end{section}

\shortorlongversion{}{
\section{An explicit isomorphism between $\QQ \lab X_n \rab$ and $ \QQ[X_n] \otimes \AA_n'.$}

Let $V$ be a finite--dimensional vector space over $\QQ$. Let $S(V)$ and $T(V)$ be its symmetric algebra and tensor algebra respectively. There exists a unique embedding $\varphi$ of $GL(V)$--modules of $S(V)$ into $T(V)$ such that 
\begin{multline*}
\varphi ( v_1 v_2 \cdots v_r ) =
\sum_{\sigma \in \S_r} v_{\sigma(1)} \otimes v_{\sigma(2)} \otimes 
\cdots \otimes v_{\sigma(r)}\\
\hbox{for all\ } r\geq 0, \; v_1, v_2, \ldots, v_r \in V. 
\end{multline*}
Its image is the subspace of the symmetric tensors. 
In the case $V=\bigoplus_{i=1}^n \QQ x_i$, we have $S(V)=\QQ[X_n]$ and $T(V)=\QQ\lab X_n \rab$. 
Then the embedding $\varphi$ and the inclusion $\AA'_n \subset \QQ\lab X_n \rab$ induce a map of $GL_n(\QQ)$--modules $\Phi: \QQ[X_n] \otimes \AA'_n \longrightarrow \QQ\lab X_n \rab$ characterized by $\Phi(f \otimes a)=\varphi(f) a$ for all $f \in \QQ[X_n]$ and all $a \in \AA'_n$. Then,
\begin{prop}[Drensky, \cite{D} Theorem 2.6]
The map $\Phi$ is a $GL_n(\QQ)$ equivariant isomorphism from $\QQ[X_n] \otimes \AA'_n$ to $\QQ \lab X_n \rab$. 
\end{prop}
Indeed, Drensky showed that given an arbitrary homogeneous basis of 
$\mathcal{G}$ of $\AA'_n$, the elements $\Phi(m \otimes g)$ for $m$ monomial 
and $g \in \mathcal{G}$, are a basis of $\AA'_n$ (\cite{D} Lemma 2.4). We refine 
Drensky's proof by considering for $\mathcal{G}$ the \emph{bracket basis} $\{P_w\}_{w \in [n]^r}$ of 
$\AA'_n$ (introduced before the proof of Proposition \ref{glniso}) and the \emph{shuffle basis} (see below) to realize $\QQ[X_n]$ 
in $\QQ\lab X_n \rab$.  We show that the elements $\Phi(m \otimes g)$ form a basis
$\QQ\lab X_n \rab$
(the \emph{hybrid basis}) that is triangularly related and expands positively in 
the \emph{bracket basis} of $\QQ \lab X_n \rab$ (Theorem \ref{refinement} below).

We follow the book of Reutenauer \cite{R} for the classical definitions
and results used in this section. The bracket basis $P_w$ has been introduced in the previous section (before the proof of Proposition \ref{glniso}).
Before presenting the hybrid basis we introduce another  
classical basis of $\QQ\lab X_n \rab$: the \emph{shuffle basis}.

\subsection*{The shuffle basis of $\QQ \lab X_n \rab$.}

Consider two monomials, $x_{i_1} x_{i_2} \cdots x_{i_r}$ and
$x_{j_1} x_{j_2} \cdots x_{j_{r'}}$ in $\QQ \lab X_n \rab$.
For a subset
\[ S = \{s_1, s_2, \ldots, s_r\}  \subseteq [r+r'], \]
and the complement subset $T = \{t_1, t_2, \ldots, t_{r'}\} = [ r+r' ]
\backslash S$, we let 
\[ x_{i_1} x_{i_2} \cdots x_{i_r}
\shuf_S x_{j_1} x_{j_2} \cdots x_{j_{r'}} := w \] be the unique monomial 
in $\QQ \lab X_n \rab$ of length $r+r'$ such
that $w_{s_1} w_{s_2} \cdots w_{s_r} = x_{i_1} x_{i_2} \cdots x_{i_r}$
and $w_{t_1} w_{t_2} \cdots w_{t_{r'}} = x_{j_1} x_{j_2} \cdots x_{j_{r'}}$.

The shuffle of any two monomials is defined as
\[ u \shuf v = \sum_{
\substack{S \subseteq [ |u|+ |v| ]
\\
|S| = |u|}
}
u \shuf_S v. \]

This shuffle of monomials is then extended to a bilinear operation on any two
elements of $\QQ\lab X_n \rab$
The shuffle product is a commutative and associative operation on
$\QQ\lab X_n \rab$.

Let $w$ be a word in $[n]^r$ and let $w = \ell_1^{i_1} \ell_2^{i_2} \cdots \ell_k^{i_k}$ be the
factorization of $w$ into decreasing products of Lyndon words $\ell_1> \ell_2 > \cdots > \ell_k$.  For a
Lyndon word $\ell=i_1 i_2 \cdots i_r$, let $S_\ell$ be the
corresponding monomial in $\QQ\lab X_n \rab$, that is $S_\ell=x_{i_1} x_{i_2} \cdots x_{i_r}$.  If $w$ is not a single Lyndon word
then define
\[ S_w = \frac{1}{i_1! i_2! \cdots i_k!}
S_{\ell_1}^{\shuf i_1} \shuf S_{\ell_2}^{\shuf i_2} \shuf \cdots \shuf S_{\ell_k}^{\shuf i_k}. \]

The set $\{ S_w \}_{w \in [n]^r}$ forms a basis
for the noncommutative polynomials of degree $r$ (\cite{R}, Corollary 5.5).

It is interesting to note
that the bracket basis $P_w$ and the shuffle basis $S_w$ are dual with respect to the scalar product where
the noncommutative monomials are self-dual.

\subsection*{The hybrid basis of $\QQ \lab X_n \rab$.} We are now ready to introduce the
  \emph{hybrid basis}.

\def\j{{a}}
Given a word $w \in [n]^r$ with a factorization into decreasing products
of Lyndon words
$w = \ell_1^{i_1} \ell_2^{i_2} \cdots \ell_k^{i_k}$, then let $\ell_{j_1},
\ell_{j_2}, \cdots, \ell_{j_r}$ be the Lyndon words of length $1$ in this decomposition
and set 
\[ M(w) = 
x_{\ell_{j_1}}^{\shuf i_{j_1}} \shuf x_{\ell_{j_2}}^{\shuf i_{j_2}} \shuf \cdots \shuf
x_{\ell_{j_r}}^{\shuf i_{j_r}} = i_{j_1}! i_{j_2}! \cdots i_{j_r}! 
S_{\ell_{j_1}^{i_{j_1}} \ell_{j_2}^{i_{j_2}} \cdots \ell_{j_r}^{i_{j_r}}}. \]
Observe that $M(w)$ is the image under the embedding $\varphi$ of the monomial 
$X(w)=\x_{\ell_{j_1}}^{i_{j_1}} \x_{\ell_{j_2}}^{i_{j_2}} \cdots 
\x_{\ell_{j_r}}^{i_{j_r}}$.
For all of the remaining Lyndon words $\ell_{\j_1}$, $\ell_{\j_2}$, \ldots,
$\ell_{\j_{k-r}}$ with length greater than $1$ we define the Lie portion of the word
to be $L(w) = P_{\ell_{\j_1}^{i_{\j_1}} \ell_{\j_2}^{i_{\j_2}} \cdots
\ell_{\j_{k-r}}^{i_{\j_{k-r}}}}$.  We will define the hybrid elements
to be $H_w := M(w) L(w)=\Phi(X(w) \otimes L(w))$. 

The result of this section is:
\begin{thm}\label{refinement}
The noncommutative polynomials $H_w$ are triangularly 
related to and expand positively in the $P_u$ basis. Precisely, for $w$ of length $r$,
\[ H_w = r!\,P_w + \hbox{terms $c_u P_u$ with $u$ lexicographically
smaller than $w$}. \] 
As a consequence, the set $\{ H_w \}_{w \in [n]^r}$ is a basis for the noncommutative polynomials
of $\QQ \lab X_n \rab$ of degree $r$.
\end{thm}
We require a few facts about Lyndon words and the lexicographic ordering
which can be found in \cite{R}.
\begin{enumerate}
\item If $u$ and $v$ are Lyndon words and $u < v$ then $uv$ is a Lyndon word. (\cite{R},
(5.1.2)) \label{fact1}
\item If $u < v$ and $u$ is not a prefix of $v$, then $ux < vy$ for all words $x,y$. (\cite{R},
Lemma 5.2.(i)) \label{fact2}
\item If $w = \ell_1 \ell_2 \cdots \ell_k$ with $\ell_1 \geq \ell_2 \geq \cdots \geq \ell_k$
then $\ell_k$ is the smallest (with respect to the $>$ order) nontrivial suffix of $w$. (\cite{R}, Lemma 7.14) \label{fact3}
\item If $\ell' < \ell$, are both Lyndon words, then $\ell' \ell < \ell \ell'$ (follows from (\ref{fact1})).
As a consequence,
for $\ell_1 \geq \ell_2 \geq \cdots \geq \ell_k$, $\ell_1 \ell_2 \cdots \ell_k \geq 
\ell_{\sigma(1)} \ell_{\sigma(2)} \cdots \ell_{\sigma(k)}$ for any permutation $\sigma \in \S_k$
with equality if and only $\ell_i = \ell_{\sigma(i)}$ for all $1\leq i \leq k$. \label{fact4}
\end{enumerate}


\begin{pf*}{Proof (of (\ref{fact4}))}  To see that (\ref{fact4}) holds  consider a weakly decreasing product of Lyndon words
$\ell_1 \ell_2 \cdots \ell_k$.  If $id \rightarrow \sigma^{(1)} \rightarrow \sigma^{(2)}
\rightarrow \cdots \rightarrow \sigma$ is a chain in the weak right
order then we have just shown that
\[ \ell_{\sigma^{(i)}(1)} \ell_{\sigma^{(i)}(1)} \cdots \ell_{\sigma^{(i)}(k)} \geq
\ell_{\sigma^{(i+1)}(1)} \ell_{\sigma^{(i+1)}(2)} \cdots \ell_{\sigma^{(i+1)}(k)} \]
with equality if and only if the two Lyndon factors which were transposed are
equal.  Therefore there exists a chain of words one greater than or equal to the
next with $\ell_1 \ell_2 \cdots \ell_k$ on one end and 
$\ell_{\sigma(1)} \ell_{\sigma(2)} \cdots \ell_{\sigma(k)}$ on the other.
\end{pf*}




Theorem \ref{refinement} will be established after the following lemma.
\begin{lem}\label{triang}
Let $w$ be a word and $\ell_1 \ell_2 \cdots \ell_r$ the decomposition
of $w$ into a decreasing product of Lyndon words.  Let
$\ell$ be a Lyndon word such that $\ell= a f_1 f_2 \cdots f_k$ with $a$ one of the variables,  each $f_i$ a Lyndon word and $f_{i} \geq f_{i+1}$ and
$f_k \geq \ell_1$.
Let $u = 
\ell_1 \cdots \ell_d \ell \ell_{d+1} \cdots \ell_r$ where $\ell_d > \ell \geq \ell_{d+1}$ or $d=0$ and $\ell \geq \ell_1$. Then
\begin{multline*}
P_\ell P_w = P_u \\+ \hbox{terms $c_v P_v$ where $v$ is lexicographically smaller than $u$ and $c_v \geq 0$}. 
\end{multline*}
\end{lem}

\begin{pf*}{Proof}
Assume that $r = 1$, and we have that either
$\ell \geq \ell_1$ and $P_\ell P_{\ell_1} = P_{\ell \ell_1}$ and we are done,
or $\ell < \ell_1$ and
\[ P_{\ell} P_{\ell_1} = P_{\ell_1} P_{\ell} + [ P_{\ell}, P_{\ell_1}]. \]
In this case $P_{\ell_1} P_{\ell} = P_{\ell_1\ell}$.   
By (\ref{fact1})
we know that
$\ell\ell_1$ is Lyndon. Moreover, $\ell \ell_1$ 
is its standard factorization (this follows from (\ref{fact3}), since the nontrivial suffixes of $\ell\ell_1$ are all suffixes of $f_1 f_2 \cdots f_k \ell_1$, which is a nonincreasing product of Lyndon words). 
Therefore $P_{\ell\ell_1} =
[P_\ell, P_{\ell_1}]$ and $P_{\ell} P_{\ell_1} = P_{\ell_1\ell} + P_{\ell\ell_1}$.
By (\ref{fact4}), $\ell \ell_1 < \ell_1\ell$ so the triangularity
relation holds.

Now for an arbitrary $r>1$ we have the same two cases.  Either $\ell \geq \ell_1$
and $P_\ell P_{\ell_1} P_{\ell_2} \cdots P_{\ell_r} = P_{\ell w}$, or $\ell < \ell_1$ and
\[ P_{\ell} P_{\ell_1} P_{\ell_2} \cdots P_{\ell_r} = P_{\ell_1} P_{\ell} P_{\ell_2} \cdots P_{\ell_r} 
+ [ P_{\ell}, P_{\ell_1}] P_{\ell_2} \cdots P_{\ell_r}. \]

Our induction hypothesis holds for $P_{\ell} P_{\ell_2} \cdots P_{\ell_r}$
since $f_k \geq \ell_1 \geq \ell_2$, hence $P_{\ell} P_{\ell_2} \cdots P_{\ell_r} =
P_{u'} + \sum_{v'<u'} c'_{v'} P_{v'}$ 
where $u' = \ell_2 \cdots \ell_d \ell \ell_{d+1} \cdots \ell_r$.
Moreover, $P_{\ell_1} P_{u'} = P_{\ell_1 u'} = P_u$ since $\ell_1 \geq \ell_2$
and $P_{\ell_1} P_{v'} = P_{\ell_1 v'}$
since $\ell_1 \geq$ any Lyndon prefix of $v'$.

Since $[P_\ell, P_{\ell_1}] = P_{\ell \ell_1}$ by (\ref{fact3}), and
$\ell_1 \geq \ell_2$, we have by the induction hypothesis that
$P_{\ell \ell_1} P_{\ell_2}  \cdots P_{\ell_r} = P_{u''} + \sum_{v''<u''} c_{v''}'' P_{v''}$
where \[u'' = \ell_2 \cdots \ell_{d'} \ell \ell_1 \ell_{d'+1} \cdots \ell_r\] with
$\ell_{d'} > \ell \ell_1 \geq \ell_{d'+1}$.  
In order to justify the induction step we also need to have that $u'' < u$.
This follows from (\ref{fact4}) since $u''$ is a permutation of the factors of
$u$ and $\ell_1 > \ell$ and $\ell$ lies to the left of $\ell_1$ in $u'$.
\end{pf*}

We are now in a position to prove Theorem \ref{refinement}.
\begin{pf*}{Proof (of Theorem \ref{refinement})}
$H_w$ is defined as the product $M(w) L(w)$ where $M(w)$ is a
a shuffle of monomials.  It expands as  $M(w) = \sum {\tilde c}_b {\bf x}_b$ with $\sum {\tilde c}_b=r!$ and where
each monomial in $M(w)$ has the same number of $x_{1}$s, $x_2$s, etc.
We fix one such monomial ${\bf x}_b$ that as follows by indexing its letters backwards: ${\bf x}_b=x_{i_{k}} x_{i_{k-1}} \cdots x_{i_1}$. 
We define inductively words $w[k]$, \ldots, $w[1]$, $w[0]$ as follows: $w[k]:=w$ and $w[j-1]$ is the word obtained from $w[j]$ by removing one of its Lyndon factors of length $1$ equal to $x_{i_j}$. Remark that $L(w[j])=L(w)$ for all $j$. Then we establish by induction on $j$ that
\begin{multline*}
x_{i_j} x_{i_{j-1}} \cdots x_{i_1} L(w)=
P_{w[j]} \\+ \hbox{terms $c_v P_v$ with $v$ lexicographically smaller than $w[j]$}
\end{multline*}
by applying Lemma \ref{triang} with $x_{i_{j}}$ for $\ell$ and  $w[j]$ for $w$.
\end{pf*}


\idiot{
\begin{section}{Appendix: table of characters and Frobenius images}
\begin{tabular}{|c|c|c|}
\hline
Module&$Gl_n(\CC)$ character&graded $S_n$-Frobenius image\\
\hline& &\\
$\QQ\lab X_n \rab \simeq T(V)$& $\frac{1}{1-(a_1+ a_2 + \cdots + a_n)} = \sum_{k \geq 0} h_{(1^k)}[A_n]$ &
$\sum_{d=0}^n \frac{q^d}{\{q;q\}_d} h_{(n-d,1^d)}[X]$\\ & & \\
\hline & & \\
$\QQ[ X_n ] \simeq S(V)$&$\prod_{i=1}^n \frac{1}{1-a_i} = \sum_{k\geq 0} h_k[A_n]$&
$h_n\left[\frac{X}{1-q}\right]$\\ & &\\
\hline& &\\
$\bigwedge(V)$& $\prod_{i=1}^n (1+a_i) = \sum_{k=0}^n e_k[A_n]$& 
$\sum_{k=0}^n q^k h_{n-k}[X] e_k[X] $\\& &$= h_n[(1-q)X] \Big|_{q \rightarrow -q}$\\
\hline&&\\
${\mathcal A}'_n$& $\frac{\prod_{i=1}^n (1-a_i)}{1-(a_1+ a_2 + \cdots + a_n)}$&
$\sum_{d=0}^{n} \frac{q^d  }{\{ q;q \}_d} h_{(n-d,1^{d})}[X(1-q)]$\\
&$= \sum_{k \geq 0} \sum_{i=0}^k (-1)^i e_{(i,1^{k-i})}[A_n]$ &\\\hline
\end{tabular}

\end{section}}
}

\end{document}